\documentclass{amsart}
\usepackage{amsfonts,amsmath,url}
\usepackage{amssymb,a4wide}
\usepackage{dsfont,color}

\usepackage[latin1]{inputenc}
\usepackage{ifthen}
\usepackage[french]{babel}

\def\Z{\mathbb Z}

\newtheorem{lemma}               {\bf Lemme}[section]

\newtheorem{theorem}     [lemma] {\bf Th\'eor\`eme}
\newtheorem{proposition} [lemma] {\bf Proposition}
\newtheorem{corollary}   [lemma] {\bf Corollaire}

\numberwithin{equation}{section}
\newcommand{\iddots}{\mathinner{%
  \mkern1mu\raise1pt\hbox{.}%
  \mkern2mu\raise4pt\hbox{.}%
  \mkern2mu\raise7pt\vbox{\kern7pt\hbox{.}}\mkern1mu}}

\newcommand{\kommentar}[1]{}

\providecommand{\sVert}[1][-1]{ \ensuremath{\mathinner{
\ifthenelse{\equal{#1}{-1}}{ 
\rvert}{}
\ifthenelse{\equal{#1}{0}}{ 
\rvert}{}
\ifthenelse{\equal{#1}{1}}{ 
\bigr\rvert}{}
\ifthenelse{\equal{#1}{2}}{ 
\Bigr\rvert}{}
\ifthenelse{\equal{#1}{3}}{ 
\biggr\rvert}{}
\ifthenelse{\equal{#1}{4}}{ 
\Biggr\rvert}{}
}} 
}

\author[Myriam Amri]{Myriam Amri}
 \address{Chair of Mathematics and Statistics, University of Leoben, A-8700 Leoben, AUSTRIA}
\email{myriam.amri@unileoben.ac.at}
\author[Lukas Spiegelhofer]{Lukas Spiegelhofer}
 \address{Institute of Discrete Mathematics and Geometry, Vienna University of Technology, 1040 Wien, Austria}
\email{lukas.spiegelhofer@tuwien.ac.at}
\author[J\"{o}rg Thuswaldner]{J\"{o}rg Thuswaldner}
 \address{Chair of Mathematics and Statistics, University of Leoben, A-8700 Leoben, AUSTRIA}
\email{joerg.thuswaldner@unileoben.ac.at}

 \address{}
\email{}
\date{\today}
\thanks{
The first and the third authors were supported by the FWF project P29910 ``Dynamics, geometry, and arithmetics of number representations''.
The second author was supported by the FWF project F5502-N26, which is a part of the Special Research Program ``Quasi Monte Carlo methods: Theory and Applications''
}

\date{\today}

\title [Sur la r\'{e}partition jointe de la repr\'{e}sentation d'Ostrowski]
{Sur la r\'{e}partition jointe de la repr\'{e}sentation d'Ostrowski dans les classes de r\'{e}sidue}
\begin{document}
\maketitle
\begin{abstract}
Pour deux entiers  $m_1, m_2\ge2$,  nous posons $\alpha_1=[0;\overline{1,m_1}]$ et $\alpha_2=[0;\overline{1,m_2}]$ et 
nous notons  respectivement par   $S_{\alpha_1}(n)$ et  $S_{\alpha_2}(n)$ les fonctions sommes des chiffres  dans les $\alpha_1$ et  $\alpha_2-$repr\'{e}sentations   d'Ostrowski de $n$.
Soient $b_1,b_2$ des entiers positifs v\'{e}rifiants $(b_1,m_1)=1$ et $(b_2,m_2)=1$, nous  obtenons une estimation avec un terme d'erreur $O(N^{1-\delta})$ pour le cardinal de l'ensemble suivant
$$\Big\{ 0\leq n<N;\ S_{\alpha_1}(n)\equiv a_1\pmod{b_1},\ S_{\alpha_2}(n)\equiv a_2\pmod{b_2}\Big\},$$
pour tous  les entiers $a_1$ et $a_2.$\\
Notre r\'{e}sultat devrait \^{e}tre   compar\'{e} \`{a} celui de B\'{e}sineau  et Kim qui ont trait\'{e} le cas des $q-$repr\'{e}sentations dans diff\'{e}rentes bases (qui sont premi\`{e}res entre elles).
\end{abstract}
\renewcommand{\abstractname}{Abstract}
\begin{abstract} 
For two distinct integers $m_1,m_2\ge2$,  we set $\alpha_1=[0;\overline{1,m_1}]$ and $\alpha_2=[0;\overline{1,m_2}]$ and
we denote  by  $S_{\alpha_1}(n)$ and $S_{\alpha_2}(n)$ respectively the  sum  of digits functions in the  Ostrowski $\alpha_1$ and $\alpha_2-$representations of $n$.
Let $b_1,b_2 $ be positive integers  satisfying   $(b_1,m_1)=1$ and $(b_2,m_2)=1$, we obtain an estimation with an error term
 $O(N^{1-\delta})$ for the cardinal of the following set
$$\Big\{ 0\leq n<N;\ S_{\alpha_1}(n)\equiv a_1\pmod{b_1},\ S_{\alpha_2}(n)\equiv a_2\pmod{b_2}\Big\},$$
for all integers $a_1$ and $a_2.$\\
Our result should be compared to that of B\'{e}sineau  and Kim who treated the case of the   $q-$representations in different bases (that are coprimes).
\end{abstract}
\section{Introduction} 
Dans tout cet article,  nous d\'{e}signons par $n$   un entier sup\'{e}rieur ou \'{e}gal \`{a} 2, $e(x)=e(2i\pi x),$
$\|x\|$ la distance du nombre r\'{e}el $x$ \`{a} l'entier le plus proche et $\lfloor x\rfloor$ d\'{e}note le plus grand entier $\leq x$ ainsi que $\{x\}$ la partie fractionnaire de $x $ et pour $m_1,m_2$ des entiers, $(m_1,m_2)$ est  le plus grand commun diviseur de ces entiers et nous notons par $|\mathcal{E}|  $ le nombre des \'{e}l\'{e}ments de l'ensemble $\mathcal{E}$  .\\
Ce papier est consacr\'{e} \`{a} l'\'{e}tude de la r\'{e}partition jointe des entiers non n\'{e}gatifs ayant deux diff\'{e}rentes $\alpha_i-$repr\'{e}sentations d'Ostrowski.
Avant de donner une d\'{e}finition exacte de notre objectif, nous voulons  examiner certains r\'{e}ultats ant\'{e}rieurs dans des sujets connexes.\\
Pour  la fonction somme des chiffres ordinaire $s_q$, Gelfond  \cite{Gelfond} a prouv\'{e}  que si $r$ est un entier $\geq2$ premier avec $(q-1)$ alors la fonction $s_{q}(n)$ est \'{e}quir\'{e}partie modulo $r.$
De plus, il a conjectur\'{e} que pour $q_1, q_2, r_1,r_2 $  des entiers $\geq 2$ o\`{u} les bases sont premi\`{e}res entre elles  avec $(r_1,q_1-1)=(r_2,q_2-1)=1$ alors il existe $\tau=\tau(q_1, q_2, r_1,r_2)>0$ tel que
\begin{equation}\label{001}
\Big| \Big\{ 0\leq n<N;\ s_{q_1}(n)\equiv d_1 \pmod{r_1},\ s_{q_2}(n)\equiv d_2\pmod{r_2}\Big\}\Big|=
\frac{N}{r_1r_2}+\mathcal O(N^{1-\tau}),
\end{equation}
pour tous  les entiers $d_1,d_2.$\\
B\'{e}sineau  \cite{Besineau} a apport\'{e} une premi\`{e}re contribution  importante mais faible \`{a} ce probl\`{e}me  avec un terme d'erreur $o(N).$ 
Vingt ans apr\'{e}s, Kim\cite{Kim} a  \'{e}tablit une r\'{e}ponse compl\`{e}te \`{a} cette conjecture o\`{u} il a m\^{e}me g\'{e}n\'{e}ralis\'{e} ce r\'{e}sultat en rempla\c{c}ant les $s_{q_i}$ par des fonctions compl\`{e}tement $q-$additives.
Lamberger et Thuswaldner \cite{Lamberger} \'{e}taient int\'{e}ress\'{e}s \`{a} \'{e}tudier la distribution de $s_G$ \'{e}tant la fonction somme des chiffres  dans une base $G$ v\'{e}rifiant une propri\'{e}t\'{e} de r\'{e}currence lin\'{e}aire, qui inclut un r\'{e}sulat analogue \`{a} celui de Gelfond ainsi qu'un r\'{e}sultat de type Erd\H{o}s-Kac.
Un exemple  connu  de cette base lin\'{e}aire r\'{e}currente est la suite de Fibonacci $F=(F_i)_{i\geq0}$ d\'{e}finit par $F_0=0, F_1=1$ et $F_{i+2}=F_{i+1}+F_{i}$ pour tout $i\geq0.$
(Notons qu'il faut commencer par le  terme $i=2$ dans la suite $F$ pour r\'{e}pondre \`{a} la d\'{e}fintion pr\'{e}c\'{e}demment donn\'{e}e). L'\'{e}tude de $S_F$ connue par la fonction de Zeckendorf est la fonction somme des chiffres dans la base de Fibonacci a \'{e}t\'{e} pr\'{e}sent\'{e}e en premier lieu dans \cite{Zeckendorf}.
Dans \cite{Lukas}, Spiegelhofer a \'{e}tudi\'{e} la relation entre la fonction de Zeckendorf $S_F$ et la fonction somme des chiffres ordinaire $s_q$ montrant que leurs valeurs sont r\'{e}parties ind\'{e}pendamment dans les classes de r\'{e}sidue donn\'{e}es (suivant Mauduit et Rivat \cite{Mauduit2}), en montrant que: Si $ \gamma \in \mathbb{R},\ \theta \in \mathbb{R}\backslash \mathbb{Z},$ alors il existe $\eta>0$ tel que
$$\sum_{ n<N}e (\theta s_{q} (n)+\gamma S_F(n))\ll N^{1-\eta}.$$
Par cons\'{e}quent, il a obtenu le r\'esultat suivant
\cite[Corollaire 5.3]{Lukas} :
Soient les entiers $q, t_2\geq2$ et  $t_1 \geq1$  tel que $(t_1,q-1)=1.$ Alors, il existe $\eta>0$ tel que pour tous les entiers $l_1, l_2$, on a
$$\Big| \Big\{ 0\leq n<N;\ s_{q}(n)\equiv l_1\pmod{t_1},\ S_F(n)\equiv l_2\pmod{t_2}\Big\}\Big|=\frac{N}{t_1t_2}+\mathcal O(N^{1-\eta}).$$
Une version non quantitative de ce r\'{e}sultat a \'{e}t\'{e} prouv\'{e}e par Coquet, Rhin et Toffin dans \cite{Rhin}. Sharma \cite{Sharma} a \'{e}tablit un r\'{e}sultat analogue \`{a} celui dans \cite{SP}  dans le cas de la fonction somme des chiffres dans la $[0;\overline{1,m}]-$repr\'{e}sentation   d'Ostrowski.  \\
Rappelons tout d'abord que si $\alpha$ est un nombre r\'{e}el, on peut le repr\'{e}senter en fraction continue simple ayant une expression de la forme
$$\alpha=a_0+\frac{1}{\displaystyle a_1+\frac{1}{\displaystyle a_2+\frac{1}{\displaystyle a_3+\cdots}}}$$
qui est d'habitude abr\'{e}g\'{e} comme $\alpha=[a_0;a_1,a_2,\dots],$
tel que $a_0=\lfloor \alpha \rfloor$ et $a_1, a_2,\ldots$ des entiers strictement positifs. La suite $(a_n)_{n\in\mathbb{N}}$ nomm\'{e}e la suite des quotients partielles  peut \^{e}tre finie ou infinie.
Si $\alpha $ est irrationnel, sa repr\'{e}sentation est infinie et unique et si $\alpha$ est rationnel, alors ils existent deux possibilit\'{e}s de repr\'{e}sentations finies. En effet, il est bien connu que
$[a_0;a_1,\ldots,a_{s-1},a_{s},1]=[a_0;a_1,\ldots,a_{s-1},a_{s}+1],$ (\`{a} voir \cite{Hardy,Perron,Shallit}).
Dans ce travail, nous sommes seulement interess\'{e}s  \`{a} $\alpha$ irrationnel entre 0 et 1. Pour une repr\'{e}sentation en fraction continue de $\alpha$ donn\'{e}e, il est possible de construire une suite de rationnels $\frac{p_i}{q_i}=[0;a_1,\ldots,a_i] $  qui est  dite la $i$-\`{e}me r\'{e}duite. Elles sont d\'{e}finies par
$p_{-1}=1,\ p_0=a_0,\ q_{-1}=0,\ q_0=1 $ et  pour tout $i\geq 0$
 $$p_{i+1}=a_{i+1}p_{i}+p_{i-1}\ \ \text{et}\ \
q_{i+1}=a_{i+1}q_{i}+q_{i-1}.$$
Il est bien connu que les suites $(p_i)_{i\geq0}$ et $ (q_i)_{i\geq0} $ v\'{e}rifient $p_{i+1}q_i-p_{i}q_{i+1}=(-1)^i$ et les $\frac{p_i}{q_i}$  convergent vers $\alpha.$ De plus, les fractions continues fournissent une suite des meilleures approximations rationnelles d'un nombre irrationnel.\\
Tout entier non n\'{e}gatif $n$ a une repr\'{e}sentation dite la {\em repr\'{e}sentation d'Ostrowski}
\begin{equation}\label{01}
n=\sum_{i\geq0}\varepsilon _i(n) q_i,
\end{equation}
sachant que $ \sum_{0\leq i<K}\varepsilon _i(n) q_i <q_K.$
Cet algorithme emm\`{e}ne \`{a} l'unicit\'{e} de la repr\'{e}sentation de la forme (\ref{01}) v\'{e}rifiant la condition dite markovienne  tel que $0\leq \varepsilon_0<a_1$ et pour tout $i\geq1$
$0\leq \varepsilon_i\leq a_{i+1}$ et quand $\varepsilon_i=a_{i+1}$ alors $\varepsilon_{i-1}=0.$
Nous r\'{e}f\'{e}rons au lecteur l'enqu\^{e}te donn\'{e}e par Berth\'e~\cite{Valerie} sur le syst\`{e}me de num\'{e}ration d'Ostrowski. Il est \`{a} noter que ce syst\`{e}me de num\'{e}ration est un outil tr\'{e}s utile pour l'\'{e}tude de la discr\'{e}pance des suites $(\alpha n) $  \`{a} voir par exemple \cite{Baxa}.
En particulier, cel\`{a} tient \`{a} la repr\'{e}sentation de Zeckendorf correspondante au cas du nombre d'Or $\alpha=\frac{\sqrt{5}-1}{2}=[0;1,1,\ldots]=[0;\overline{1}],$ o\`{u} la suite
$(q_i)_{i\geq0} $ correspond aux nombres de Fibonacci \cite{Vorobiev}. L'unicit\'{e} de la condition $\varepsilon_i=0$ si $\varepsilon_{i+1}=a_{i+1}$ signifie qu'on ne peut pas avoir deux termes cons\'{e}cutifs dans la repr\'{e}sentation de Zeckendorf correspondante \cite{Zeckendorf}.
Pour plus de d\'{e}tails, voir par exemple \cite{Allouche,Shallit}.\\
Le but de ce papier est de prouver un r\'{e}sultat analogue \`{a} celui de Kim \cite{Kim}
pour le syst\`{e}me de num\'{e}ration d'Ostrowski.\\
Donnons un entier $m\geq2$ et un  irrationnel $\alpha\in (0,1)$ tel que sa fraction continue est de la forme  $\alpha=[0;\overline{1,m}]$. Alors $q_0=q_0(m)=q_1=q_1(m)=1 $ et
\begin{equation}\label{def:qi}
q_i=q_i(m)=
 \left\{
\begin{array}{lll}
mq_{i-1}+q_{i-2},  &si& i\ \text{est pair} \\
q_{i-1}+q_{i-2},  &si& i\ \text{est impair}\\
\end{array} \right.
\end{equation}
et
\begin{equation}\label{eq:phidef}
\varphi=\varphi(m)=\frac{m+2+\sqrt{m^2+4m}}{2}.
\end{equation}

Apr\`{e}s cette bri\`{e}ve revue, nous  pr\'{e}sentons  notre Th\'{e}or\`{e}me principal.
\begin{theorem}\label{th}

Soient $m_1,m_2\geq 2$ deux entiers distincts,  $\alpha_1=[0;\overline{1,m_1}]$  et $\alpha_2=[0;\overline{1,m_2}]$. Soient $\vartheta\in \mathbb{R} $ et $m_2\beta\in \mathbb{R}\setminus\mathbb{Z}$. Alors, il existe $\delta>0$
 \begin{equation}\label{002}
 \sum_{ n<N}e (\vartheta S_{\alpha_1}(n)+\beta S_{\alpha_2}(n))\ll N^{1-\delta}.
 \end{equation}
\end{theorem}
Comme simple cons\'{e}quence, nous r\'{e}sultons le corollaire suivant
\begin{corollary}\label{cor}
Soient $m_1,m_2\geq 2$ deux entiers distincts, $\alpha_1=[0;\overline{1,m_1}]$, $\alpha_2=[0;\overline{1,m_2}]$
et $b_1,b_2 $ des entiers positifs tel que $(b_1,m_1)=1$ ou $(b_2,m_2)=1$.  Alors, il existe $\delta>0$ tel que
\begin{equation*}
\Big| \Big\{ 0\leq n<N;\ S_{\alpha_1}(n)\equiv a_1\pmod{b_1},\ S_{\alpha_2}(n)\equiv a_2\pmod{b_2}\Big\}\Big|=
\frac{N}{b_1b_2}+\mathcal O(N^{1-\delta}),
\end{equation*}
pour tous  les entiers $a_1$ et $a_2.$
\end{corollary}
Le plan de ce papier est organis\'{e} comme suit: Dans la deuxi\`{e}me section, nous  introduisons quelques notations et nous  pr\'{e}sentons  une s\'{e}rie de r\'{e}sultats auxili\`{e}res qui  constituent  des  ingr\'{e}dients essentiels dans la r\'{e}solution de notre objectif  puis  nous pr\'{e}c\'{e}dons \`{a} la preuve du Th\'{e}or\`{e}me \ref{th} ainsi que celle du Corollaire \ref{cor}.
La troisi\`{e}me section contient les preuves de ces r\'{e}sultats.

\section{Pr\'{e}liminaires}\label{sec:prelim}
Le lemme suivant est un r\'{e}sultat \'{e}l\'{e}mentaire sur les sommes d'exponentielles.
\begin{lemma}[\cite{Korobov}]\label{Korobov}
Soient $x\in \mathbb{R}$ et $N, R\geq 0.$ Alors, on a
\begin{equation}\label{02}
  \sum\limits_{r<R} \Big(R-|r|\Big)e(rx)=\Big|\sum\limits_{r<R} e(rx)\Big|^2  .
\end{equation}
\end{lemma}
\begin{lemma}[{\cite[Corollaire\footnote{Il est \`{a} noter qu'il a une faute de frappe concernant le sens de l'in\'{e}galit\'{e} dans \cite[Corollaire]{Schmidt}}]{Schmidt}\label{lem2.2}}]
Soient $\alpha_1,\ldots,\alpha_n$ des nombres r\'{e}els alg\'{e}briques tel que $1,\alpha_1,\ldots,\alpha_n$ sont lin\'{e}airement ind\'{e}pendants  sur l'ensemble des rationnels $\mathbb{Q}$. Alors, pour tout $\varepsilon>0,$
il existe une constante $c>0$, pour tout $h_1,\ldots,h_n,p$ avec $(h_1,\ldots, h_n) \not=(0,\ldots,0),$ tel que
\begin{equation}\label{lemsch}
\Big\|h_1\alpha_1+\ldots+ h_n\alpha_n+p\Big\|\geq c h^{-n-\varepsilon},
\end{equation}
o\`{u}  $h=\max(|h_1|,\ldots,|h_n|)>0.$
\end{lemma}
Nous notons maintenant une estimation analogue \`{a} \cite[Lemma 5.8]{SP}, prouv\'{e} en utilisant la discr\'{e}pance de la suite $(n\varphi)$ ayant une quotient partielle born\'{e}e.

\begin{lemma}[{\cite[Lemma 5.8]{SP}}]\label{lem2}
Soit $I$ un interval fini inclus dans $\mathbb{Z}$. On suppose que $K$ et $t$ sont des nombres r\'{e}els et $K\geq 1$. Alors, nous avons
$$\sum\limits_{h\in I} \min \Big\{K, \frac{1}{\|t+h\varphi \|^2} \Big\}
\ll \sqrt{K} |I|+K\log|I|.$$
\end{lemma}
La preuve du Th\'{e}or\`{e}me \ref{th} n\'{e}cessite
la  version suivante de l'in\'{e}galit\'{e} de Weyl-Van der Corput.
\begin{lemma}[{\cite[Lemma 2.5]{Graham}}] \label{l4}
Soient les nombres complexes  $a_0,a_1,\ldots,a_{N-1} $ . Alors, pour tout entier $R\geq1$, nous avons
$$ \Big|\sum\limits_{0\leq n<N}a_n \Big|^2\leq \frac{N-1+R}{R}\sum\limits_{0\leq |r|<R}\Big(1-\frac{|r|}{R}\Big)\sum_{0\leq n<N\atop{0\leq n+r<N}} a_{n+r} \overline{a}_n.$$
\end{lemma}

Soit un entier $k\geq2$, nous notons par  $t(n,k)$ la troncature de la repr\'{e}sentation d'Ostrowski de $n$  c'est \`{a}  dire
$$t(n,k)=\sum_{0\leq i<k}\varepsilon _i(n) q_i$$
et   par
$S_{\alpha,k}$ la fonction somme des chiffres jusqu'au terme $k$ c'est \`{a} dire
$$S_{\alpha,k}(n)=\sum_{0\leq i<k}\varepsilon _i(n).$$
Le lemme suivant concerne les fonctions sommes des chiffres tronqu\'{e}es pour prouver l'id\'{e}e de l'addition d'un entier  $r$ \`{a}  $n$ devant changer les digits dans des positions tr\`{e}s basses dans la plus part des cas.

\begin{lemma}[{\cite[Lemma 2.4]{Sharma}}] \label{lem1}
Soient $N, k$ et $r$ des entiers non n\'{e}gatifs tel que $k\geq 2$. Alors, nous avons
$$\Big|\Big\{ n<N;\; S_{\alpha}(n+r)-S_{\alpha}(n)\neq S_{\alpha,k}(n+r)-S_{\alpha,k}(n) \Big\}\Big|\leq \frac{Nr}{ q_{k-1}}.$$
\end{lemma}

\begin{lemma}[{\cite[Lemma 3.5]{Sharma}}] \label{lem01}
Soient   $\gamma\in \mathbb{R}$ avec $\|m\gamma\|\neq0 $. Alors, ils existent $c$ et $\eta >0$ sachant que pour tout  $\theta\in \mathbb{R}$ et $k\geq2,$ nous avons
$$\Big|\frac{1}{{q}_k}\sum_{0\leq u< q_{k}}e\Big(\gamma S_{\alpha}(u)+ \theta u\Big)\Big|\leq  \eta e^{-ck}.$$
\end{lemma}
La clef de la d\'{e}monstration du Th\`{e}or\`{e}me~\ref{th} r\'{e}side dans l'obtention de l'estimation de l'exponentielle de la fonction somme des chiffres tronqu\'{e}e suivante.

\begin{proposition}[{\cite[Lemma 3.4]{Sharma}}]\label{prop1}
Soient $n$ un entier non n\'{e}gatif, $h\in \mathbb{Z}$ et $\vartheta \in \mathbb{R}.$
Pour H et k des entiers  avec
$k\geq 2$ et $H\geq1$, nous d\'{e}finissons
$$M^{(1)}_k(h,\vartheta)=\sum\limits_{0\leq u< q_{k-1}} e \Big(\vartheta S_{\alpha} (u) -(-1)^k hu\varphi\Big),$$
$$M^{(2)}_k(h,\vartheta)=\sum\limits_{q_{k-1}\leq u< q_{k}} e \Big(\vartheta S_{\alpha} (u) -(-1)^k hu\varphi \Big). $$
Pour $|h|\leq H$, ils existent des nombres complexes $b_{H} ^{(1)}(h) , b_{H} ^{(2)}(h),
c_{H}^{(1)}(h)$ et $c_{H}^{(2)}(h)$ tel que pour
\begin{equation*}
b_{H} ^{(1)}(0)=
 \left\{
\begin{array}{lll}
\frac{2-m+\sqrt{d}}{2\varphi^{k_0}},  &si&k=2k_0, k_0\in \mathbb{Z}\\
\frac{1}{\varphi ^{k_0}},  &si&k=2k_0+1, k_0\in \mathbb{Z}\\
\end{array} \right.
\end{equation*}
\begin{equation*}
b_{H} ^{(2)}(0)=
 \left\{
\begin{array}{lll}
\frac{1}{\varphi ^{k_0}},  &si&k=2k_0, k_0\in \mathbb{Z}\\
\frac{-m+\sqrt{d}}{2\varphi^{k_0}},  &si&k=2k_0+1, k_0\in \mathbb{Z}\\
\end{array} \right.
\end{equation*}
et pour tout $i=1,2$
\begin{equation}\label{10}
|b_{H} ^{(i)}(h)|\leq\min \Big(b_{H} ^{(i)}(0), \frac{1}{|h|}\Big),\quad (\text{si}\  h\neq 0) \end{equation}
$$|c_{H}^{(i)} (h)|\leq 2,$$
on a
\begin{eqnarray}\label{p2}
\begin{split}
e (\vartheta S_{\alpha,k}(n))=& \sum\limits_{|h|\leq H} e \Big((-1)^k hn\varphi\Big)b_ {H}^{(1)} (h)M^{(1)}_k(h,\vartheta)+\sum\limits_{|h|\leq H} e \Big((-1)^k hn\varphi \Big)b_ {H}^{(2)} (h)M^{(2)}_k(h,\vartheta) \nonumber\\
&+\mathcal O\Big(\frac{1}{H}\sum\limits_{|h|\leq H}c_ {H}^{(1)} (h) e \Big( (-1)^k hn\varphi \Big) \sum\limits_{0\leq u< q_{k-1}} e \Big( -(-1)^k hu\varphi \Big)\Big)\nonumber  \\
 &+\mathcal O\Big(\frac{1}{H}\sum\limits_{|h|\leq H}c_ {H}^{(2)} (h) e \Big( (-1)^k hn\varphi \Big) \sum\limits_{q_{k-1}\leq u< q_{k}} e \Big( -(-1)^k hu\varphi\Big) \Big),
\end{split}
\end{eqnarray}
sachant que les termes d'erreurs  sont des nombres r\'{e}els non n\'{e}gatifs et les constantes implicites sont absolues.
\end{proposition}
Nous  voulons d\'{e}velopper quelques changements sur  $S_{\alpha, k}$ en des s\'{e}ries de Fourier ordinaires. Nous posons
$V=\{ n\in \mathbb{N},\ \varepsilon_j(n)=0\ \hbox{pour}\ j<k\}.$
Soit la suite $(n_v)_{v\geq0}$ obtenue en ordonnant les \'{e}l\'{e}ments de $V$ dans l'ordre croissant. Il est clair que  $n_0=0$ et $Q(v)=n_v-n_{v-1}\in
\{q_{k-1}, q_{k}\},$  pour tout $v\geq 1.$
Nous consid\'{e}rons la transform\'{e}e de Fourier discr\`{e}te suivante
$$L_v(l)=\frac{1}{Q(v)}\sum\limits_{u<Q(v)}  e \Big(\vartheta S_{\alpha,k} (u+n_{v-1})-  \frac{lu}{Q(v)}\Big).$$
Il suit que
\begin{equation}\label{07}
e (\vartheta S_{\alpha,k}(n+n_{v-1}))= \sum\limits_{l<Q(v)} L_v(l) e\Big(\frac{ln}{Q(v)}\Big).
\end{equation}
Notons que la somme  (\ref{07}) est  \`{a} priori vrai pour tout   $0\leq n<Q(v). $ De plus puisque  $S_{\alpha,k}(n+n_{v-1})=S_{\alpha,k}(n+n_{v-1}+Q(v))$ pour tout  $0\leq n<q_{k-1}$.
On constate que  (\ref{07}) est aussi vraie pour  $0\leq n<Q(v)+q_{k-1}$.\\

\section{Preuve du Th\'{e}or\`{e}me~\ref{th}}

Dans la suite, nous utilisons les abr\'eviations
$\varphi_1=\varphi(m_1)$ et $\varphi_2=\varphi(m_2)$ o\`{u} $\varphi(m)$ est d\'{e}finit dans \eqref{eq:phidef}. De plus, nous notons par  $q_{1,k}=q_{k}(m_1)$ et $q_{2,k}=q_{k}(m_2)$ o\`{u} $q_{k}$ est d\'{e}finit dans  \eqref{def:qi}. Il est \`{a} indiquer  aussi que la propri\'{e}t\'{e} suivante
\[
\varphi_i^k \ll  q_{i,k} \ll \varphi_i^k
\]
est d\'{e}tenue  \`{a } cause de la relation de r\'{e}currence de $q_{i,k}$ ($i\in\{1,2\}$).  Celle-ci  sera ult\'{e}rieurement utilis\'{e}e dans le cadre de notre preuve.\\
Nous commen\c{c}ons par l'application du Lemme~\ref{l4}  sur la somme d\'{e}sir\'{e}e. Alors, nous  obtenons
\begin{equation}\label{1}
\begin{split}
\Big|\sum_{ n<N}e (\vartheta S_{\alpha_1}(n)+&\beta S_{\alpha_2}(n))\Big|^2 \ll\frac{N-1+R}{R}\sum\limits_{0\leq |r|<R}\Big(1-\frac{|r|}{R}\Big) \\
 &\sum_{0\leq n\leq N\atop{0\leq n+r\leq N}} e \Big(\vartheta \big(S_{\alpha_1}(n+r)- S_{\alpha_1}(n)\big)\Big)\  e \Big(\beta \big(S_{\alpha_2}(n+r)- S_{\alpha_2}(n)\big)\Big).
 \end{split}
\end{equation}
En premier lieu, le terme d'erreur d'ordre $\mathcal O(NR)$ d\'{e}coule de la n\'{e}gligence
de la condition  $0\leq n+r\leq N$ et du remplacement du terme $N-1+R $ par $N$.
Afin d'estimer la  somme \`{a} droite, nous aurons besoin  des fonctions tronqu\'{e}es.
Alors, un  appel au Lemme \ref{lem1} nous permet de
remplacer  $S_{\alpha_i}$ par $S_{\alpha_i,k_i}$ pour $i=1,2$ dans (\ref{1}) en engendrant les termes d'erreur d'ordre $\mathcal O\Big(\frac{N^2R}{q_{1,k_1-1}}\Big)$ et $\mathcal O\Big(\frac{N^2R}{q_{1,k_2-1}}\Big).$
En assemblant tous ces termes d'erreur, nous parvenons \`{a} \'{e}crire
\begin{equation}\label{eqp}
\begin{split}
            \Big|\sum_{ n<N}e (\vartheta S_{\alpha_1}(n)+\beta S_{\alpha_2}(n))\Big|^2\ll&S
            +\mathcal O(NR)+\mathcal O\Big(\frac{N^2R}{q_{1,k_1-1}}+\frac{N^2R}{q_{2,k_2-1}}\Big),
\end{split}
\end{equation}
o\`{u}
$$S=\frac{N}{R}\sum\limits_{0\leq |r|<R}\Big(1-\frac{|r|}{R}\Big)
            \sum_{n\leq N}  e \Big(\vartheta (S_{\alpha_1,k_1}(n+r)- S_{\alpha_1,k_1}(n))\Big)
             e \Big(\beta (S_{\alpha_2,k_2}(n+r)- S_{\alpha_2,k_2}(n))\Big).$$
Il suit qu'en rempla\c{c}ant les facteurs  d'exponentielles
\[
E_2 = e (\beta S_{\alpha_2,k_2}(n+r))e (- \beta S_{\alpha_2,k_2}(n))
\]
par leurs expressions donn\'{e}es dans la proposition pr\'{e}c\'{e}dante, nous obtenons seize sommes de produits de  termes principaux et de termes d'erreurs dont nous  classifions en trois classes diff\'{e}rentes:
produit de deux termes principaux, produit de deux termes d'erreurs et produit d'un terme principal et un terme d'erreur. Nous consid\'{e}rons ces diff\'{e}rents cas d'une fa\c{c}on s\'{e}par\'{e}e. Par cons\'{e}quent, nous d\'{e}visons la somme  $S$   en trois parties
\[
S=S_1+S_2+S_3.
\]

\subsection{Cas 1: Si les deux facteurs contribu\'es par $E_2$ sont des termes d'erreurs} \label{sec:Cas1}
Alors nous avons la contribution suivante
\begin{eqnarray}\label{eq:S0_1}
\begin{split}
S_1=&\frac{N}{R}\sum\limits_{0\leq |r|<R}\Big(1-\frac{|r|}{R}\Big)\sum_{n\leq N} e \Big(\vartheta (S_{\alpha_1,k_1}(n+r)- S_{\alpha_1,k_1}(n))\Big) \\
&  \cdot \mathcal O\Big(\frac{1}{H}\sum\limits_{|h_1|\leq H}c_ {H}^{(j_1)} (h_1) \sum\limits_{u} e ( (-1)^{k_2} h_1 \varphi_2(n+r-u))\Big)\\
& \cdot \mathcal O\Big(\frac{1}{H}\sum\limits_{|h_2|\leq H}c_ {H}^{(j_2)} (h_2) \sum\limits_{u} e \Big( (-1)^{k_2} h_2 \varphi_2(n-u)\Big)\Big),
\end{split}
\end{eqnarray}
pour tout  $0\leq u<q_{2,k_2-1}$ ou $q_{2,k_2-1}\leq u<q_{2,k_2}$ et $j_1,j_2\in\{1,2\}$.
Ensuite, en prenant $e (\beta S_{\alpha_1,k_1}(n+r))$ ainsi que $e (-\beta S_{\alpha_1,k_1}(n))$ et les rempla\c{c}ant chacune par leurs formes en (\ref{07}) et les ins\'{e}rant dans la somme sur $n$.
Pour cel\`a, nous supposons que $R<\frac{q_{1,k_1}}{2}$ et nous choisissons $v_1$ avec $n_{v_1-1}\le N< n_{v_1}$. Si nous  changeons l'intervalle de la sommation sur $n$ dans \eqref{eq:S0_1}  \`a $0\le n-R \le n_{v_1-1}$ , \c{c}a nous co\^ute un terme d'erreur $\mathcal O(Nq_{1,k_1}q_{2,k_2}^2)$et si  on bloque cette nouvelle somme sur $n$ en prenant en consid\'eration la suite $(n_v)$, on gagne
\begin{eqnarray*}
\begin{split}
S_1=&\frac{N}{R}\sum\limits_{0\leq |r|<R}\Big(1-\frac{|r|}{R}\Big)\sum_{v=0}^{v_1-1}\sum_{n=0}^{Q(v)-1} e \Big(\vartheta (S_{\alpha_1,k_1}(n+r+R+n_v)- S_{\alpha_1,k_1}(n+R+n_v))\Big) \\
&  \cdot \mathcal O\Big(\frac{1}{H}\sum\limits_{|h_1|\leq H}c_ {H}^{(j_1)} (h_1) \sum\limits_{u} e ( (-1)^{k_2} h_1 \varphi_2(n+R+n_v+r-u))\Big)\\
& \cdot \mathcal O\Big(\frac{1}{H}\sum\limits_{|h_2|\leq H}c_ {H}^{(j_2)} (h_2) \sum\limits_{u} e \Big( (-1)^{k_2} h_2 \varphi_2(n+R+n_v-u)\Big)\Big)+ \mathcal O(Nq_{1,k_1}q_{2,k_2}^2).
\end{split}
\end{eqnarray*}
Maintenant, on applique (\ref{07}) \`a $e (\beta S_{\alpha_1,k_1}(n+r))$ et $e (-\beta S_{\alpha_1,k_1}(n))$ pour parvenir \`{a}
\begin{eqnarray*}
\begin{split}
S_1=\frac{N}{R}\sum\limits_{0\leq |r|<R}&\Big(1-\frac{|r|}{R}\Big)\\
&\sum_{v=0}^{v_1-1}\sum_{{0\leq l_1,l_2<Q(v)}} L_v(l_1)\overline{L_v(-l_2)} e\Big(\frac{l_1r}{Q(v)}\Big)\sum_{n=0}^{Q(v)-1} e\Big(\frac{(l_1+l_2)(n+R+n_v)}{Q(v)}\Big)\\
& \cdot \mathcal O\Big(\frac{1}{H}\sum\limits_{|h_1|\leq H}c_ {H}^{(j_1)} (h_1) \sum\limits_{u} e ( (-1)^{k_2} h_1 (n+R+n_v+r-u)\varphi_2)\Big)\\
& \cdot \mathcal O\Big(\frac{1}{H}\sum\limits_{|h_2|\leq H}c_ {H}^{(j_2)} (h_2) \sum\limits_{u} e \Big( (-1)^{k_2} h_2 (n+R+n_v-u)\varphi_2\Big)\Big)+ \mathcal O(Nq_{1,k_1}q_{2,k_2}^2).
\end{split}
\end{eqnarray*}
Par suite, il suffit  de majorer trivialement sur la somme sur  $r$  ainsi que le premier terme d'erreur par $q_{2,k_2}$ (en tenant compte du fait que les expressions dans les  termes d'erreur  sont des nombres r\'{e}els non n\'{e}gatifs et les constantes implicites sont absolues) et consid\'{e}rer l'identit\'{e} de Parseval, pour obtenir la majoration suivante
\begin{eqnarray*}
S_1   &\ll&\frac{N}{H}q_{2,k_2}\sum_{v=0}^{v_1-1} \sum_{l=0}^{Q(v)-1} |L_v(l)|^2\sum\limits_{|h_2|\leq H} \Big|\sum_{n=0}^{Q(v)-1}
e \Big( h_2n  \varphi_2\Big)\Big|  \Big|\sum\limits_{u}e \Big( h_2u_2\varphi_2  \Big)\Big| + Nq_{1,k_1}q_{2,k_2}^2\nonumber\\
   &\ll&\frac{N}{H}q_{2,k_2}\sum_{v=0}^{v_1-1}
\sum\limits_{|h_2|\leq H} \min \Big(\|h_2\varphi_2\|^{-1},|Q(v)|\Big)
\min\Big(\|h_2\varphi_2\|^{-1},q_{2,k_2}\Big)+ Nq_{1,k_1}q_{2,k_2}^2.\nonumber
\end{eqnarray*}
Par ailleurs, la  majoration de $Q(v)$ par $q_{1,k_1}$ et  l'appel au lemme \ref{lem2}   nous permettent  d'\'{e}crire
\begin{eqnarray*}
S_1   &\ll&\frac{N}{H}q_{2,k_2}v_1
\sum\limits_{|h_2|\leq H} \min \Big(\Big\|h_2\varphi_2\Big\|^{-2},q_{1,k_1}q_{2,k_2}\Big)+ Nq_{1,k_1}q_{2,k_2}^2
\nonumber\\
&\ll&\frac{N}{H}q_{2,k_2}v_1\Big(\sqrt{q_{1,k_1}q_{2,k_2}}H+q_{1,k_1}q_{2,k_2}\log H \Big)+ Nq_{1,k_1}q_{2,k_2}^2.
\nonumber
\end{eqnarray*}
Il suit, qu'en utilisant le fait que  $v_1\leq c\frac{N}{q_{1,k_1}}$,
\begin{equation}\label{eq:S1final}
S_1\ll N^2 q^{-1/2}_{1,k_1}q^{3/2}_{2,k_2}+N^2q^{2}_{2,k_2}\frac{\log H}{H}+ Nq_{1,k_1}q_{2,k_2}^2.
\end{equation}

\subsection{Cas 2 :
Si les deux facteurs contribu\'es par $E_2$ sont des termes principaux} Nous avons, pour tout  $i,j\in\{1,2\}$,
 \begin{eqnarray}
\begin{split}
   &S_2=\frac{N}{R}\sum\limits_{0\leq |r|<R}\Big(1-\frac{|r|}{R}\Big)
\sum_{n\leq N} e \Big(\vartheta (S_{\alpha_1,k_1}(n+r)- S_{\alpha_1,k_1}(n))\Big)\nonumber\\
&\;\; \cdot \sum\limits_{|h_1|\leq H} b_ {H}^{(i)} (h_1)M^{(i)}_{k_2}(h_1,\beta)
e \Big((-1)^{k_2} h_1(n+r)\varphi_2\Big) \cdot \sum\limits_{|h_2|\leq H} \overline{b_ {H}^{(j)} (-h_2)} \overline{M^{(j)}_{k_2}(-h_2,\beta)}
e \Big((-1)^{k_2} h_2n\varphi_2\Big).\nonumber
\end{split}
\end{eqnarray}
Afin d'estimer $S_2$, nous devons  discuter selon la somme des termes  $h_1+h_2$.  En effet, nous divisons la somme en question en tant que   $$S_2= S_{2,0}+S_{2,1},$$ o\`{u}  
$S_{2,0}$ et $S_{2,1}$  contiennent respectivement  les sommes correspondantes au choix  $h_1+h_2=0$ et
  $h_1+h_2\neq 0.$
\subsubsection{Estimation de $S_{2,0}$ o\`{u} $h_1+h_2= 0$}
Vu que $h_1=-h_2$, nous avons
\begin{eqnarray}\label{5}
\begin{split}
S_{2,0}=\frac{N}{R}&\sum\limits_{0\leq |r|<R}\Big(1-\frac{|r|}{R}\Big)
\sum_{n\leq N} e \Big(\vartheta (S_{\alpha_1,k_1}(n+r)- S_{\alpha_1,k_1}(n))\Big)\nonumber\\
&\cdot\sum\limits_{|h_1|\leq H} b_ {H}^{(i)} (h)M^{(i)}_{k_2}(h_1,\beta)
e \Big((-1)^{k_2} h_1r\varphi_2\Big)
\overline{b_ {H}^{(j)} (h_1)} \overline{M^{(j)}_{k_2}(h_1,\beta)}.\nonumber
\end{split}
\end{eqnarray}
Nous bloquons la somme sur  $n$ de la m\^{e}me mani\`{e}re que dans la   Section~\ref{sec:Cas1} et  nous ins\'{e}rons l'expression (\ref{07}) dans la somme sur $n$. Alors, nous acqu\'{e}rons
 \begin{eqnarray*}\label{6}
\begin{split}
S_{2,0}=&\frac{N}{R}\sum\limits_{0\leq |r|<R}\Big(1-\frac{|r|}{R}\Big)
\sum_{v=0}^{v_1-1} \\
&
 \sum_{{0\leq l_1,l_2<Q(v)}} L_v(l_1)\overline{L_v(-l_2)} e\Big(\frac{l_1r}{Q(v)}\Big)\sum_{n=0}^{Q(v)-1} e\Big(\frac{(l_1+l_2)(n+R+n_v)}{Q(v)}\Big)
\\
&
\cdot \sum\limits_{|h_1|\leq H} b_ {H}^{(i)} (h)M^{(i)}_{k_2}(h_1,\beta)
e \Big((-1)^{k_2} h_1r\varphi_2\Big)
\overline{b_ {H}^{(j)}(-h_1)}\overline{M^{(j)}_{k_2}(-h_1,\beta)}+\mathcal O(NHq_{1,k_1})\\
=&\frac{N}{R^2}\sum_{v=0}^{v_1-1}\sum\limits_{{|h_1|\leq H}\atop{0\leq l_1,l_2<Q(v)}}L_v(l_1)\overline{L_v(-l_2)} b_ {H}^{(i)}(h_1) \overline{b_{H}^{(j)}(-h_1)}
M^{(i)}_{k_2}(h_1,\beta)\overline{M^{(j)}_{k_2}(-h_1,\beta)}\\
&
\cdot\sum\limits_{0\leq |r|<R}\Big(R-|r|\Big)
e\Big(r\Big(\frac{l_2}{Q(v)}+ (-1)^{k_2} h_1 \varphi_2\Big)\Big)\sum_{n=0}^{Q(v)-1} e\Big(\frac{l_1+l_2}{Q(v)}(n+R+n_v)\Big)\\
&+\mathcal O(NHq_{1,k_1}).
\end{split}
\end{eqnarray*}
Si $\ell_1+\ell_2\not\equiv 0\pmod{Q(v)}$, la somme int\'{e}rieure sur  $n$ est disparue.   Donc, nous pouvons supposer que  $\ell_1+\ell_2\equiv 0\pmod{Q(v)}$.
Lemme~\ref{Korobov}, l'estimation sur  $|b_ {H}^{(i)} (h_1)|$ dans  \eqref{10}  ainsi que l'estimation triviale  $|M^{(i)}_{k_2}(h_1,\beta)|\le q_{2,k_2}$ produisent par cons\'{e}quent
\[
\begin{split}
S_{2,0} \ll
\frac{N^2}{R^2} \sup_{h\in \mathbb{Z}}\Big|b_ {H}^{(i)}(h) M^{(i)}_{k_2}(h,\vartheta)\Big|
 \sum_{l=0}^{Q(v)-1}& |L_v(l)|^2 \sum\limits_{|h_1|\leq H} q_{2,k_2}
 \min\Big(\frac{1}{|h_1|},\varphi^{-k_2+1}_{2}\Big)\\
 \cdot&
 \min \Big( R^2, \Big\|\frac{l}{Q(v)}+ (-1)^{k_2} h_1\varphi_2\Big\|^{-2}\Big)+NHq_{1,k_1}.
\end{split}
\]
Maitenant, un appel \`{a} l'identit\'{e} de Parseval  et au  Lemme~\ref{lem01} donnent
\begin{equation}\label{eq:S20final}
\begin{split}
S_{2,0}\ll&
\frac{N^2}{R^2} e^{-c_2k_2}q_{2,k_2}
 \sup_{l\in\Z}\bigg(
\sum\limits_{|h_1|\leq H}
\min\Big(\frac{1}{|h_1|},\varphi^{-k_2+1}_{2}\Big)
 \min \Big( R^2, \Big\|\frac{l}{Q(v)}+ |h_1|\varphi_2\Big\|^{-2}\Big) \bigg)\\
 &+NHq_{1,k_1}
\\
\ll&\frac{N^2}{R^2} e^{-c_2k_2}q_{2,k_2}
\sum\limits_{{s\leq \frac{H}{q_{2,k_2}}}}\min\Big(\frac{1}{sq_{2,k_2}},\frac{1}{q_{2,k_2}}\Big)\sup_{t\in \mathbb{R}} \bigg( \sum\limits_{{h\leq q_{2,k_2}}}
\min \Big( R^2, \Big\|t+(sq_{2,k_2}+h)\varphi_2\Big\|^{-2}\Big) \bigg)\\
&+NHq_{1,k_1}
 \\
\ll&N^2 e^{-c_2k_2}\log H\Big(\frac{q_{2,k_2}}{R} + \log q_{2,k_2} \Big) +NHq_{1,k_1},
\end{split}
\end{equation}
sachant que la derni\`{e}re in\'{e}galit\'{e} r\'{e}sulte du Lemme \ref{lem2}.\\
\`{A} vraie dire, l'utilisation du Lemme~\ref{lem01} dans l'estimation d\'{e}sir\'{e}e interviendra ult\'{e}rieurement dans les conditions au Corollaire~\ref{cor}.
\subsubsection{Estimation de $S_{2,1}$ o\`{u} $h_1+h_2\neq 0$}\label{sec:h1h1not0}
\[
 \begin{split}
&S_{2,1}=\frac{N}{R}\sum\limits_{0\leq |r|<R}\Big(1-\frac{|r|}{R}\Big)
\sum_{n\leq N} e \Big(\vartheta (S_{\alpha_1,k_1}(n+r)- S_{\alpha_1,k_1}(n))\Big)\\
&
\quad \cdot  \sum\limits_{{|h_1|, |h_2|\leq H}\atop{h_1+h_2\neq0}} b_ {H}^{(i)} (h_1)M^{(i)}_{k_2}(h_1,\beta) e \Big((-1)^{k_2} h_1(n+r)\varphi_2\Big) \overline{b_ {H}^{(j)} (-h_2)} \overline{M^{(j)}_{k_2}(-h_2,\beta)}
e \Big((-1)^{k_2} h_2n\varphi_2\Big)
\end{split}
\]
En rempla\c{c}ant les facteurs  d'exponentielles
\[
E_1 = e (\beta S_{\alpha_1,k_1}(n+r))e (- \beta S_{\alpha_1,k_1}(n))
\]
 par leurs expressions dans la proposition \ref{prop1}, ils nous viennent une autre fois pour ce sous-cas seize sommes de produits de  termes principaux et de termes d'erreurs dont nous  distinguons trois cas
 \[
 S_{2,1} = S_{2,1,1} + S_{2,12} + S_{2,1,3}
 \]
comme suit.

\medskip

{\parindent 0pt \it \ref{sec:h1h1not0}.1. Si les deux facteurs contribu\'es par $E_1$ sont des termes d'erreurs.}
Dans ces cas, nous avons
 \begin{equation*}
\begin{split}
   S_{2,1,1}\ll&\frac{N}{R}\sum\limits_{0\leq |r|<R}\Big(1-\frac{|r|}{R}\Big)\nonumber\\
&\sum_{n\leq N} \sum\limits_{{|h_1|, |h_2|\leq H}\atop{h_1+h_2\neq0}} b_ {H}^{(i_1)} (h_1)M^{(i_1)}_{k_2}(h_1,\beta)e \Big((-1)^{k_2} h_1(n+r)\varphi_2\Big)\\
&
\qquad\qquad\quad\;\;\; \cdot  \overline{b_ {H}^{(i_2)} (-h_2)} \overline{M^{(i_2)}_{k_2}(-h_2,\beta)}e \Big((-1)^{k_2} h_2n\varphi_2\Big)\\
& \cdot \mathcal O\Big(\frac{1}{H}\sum\limits_{|h_3|\leq H}c_ {H}^{(j_1)} (h_3) \sum\limits_{u} e ((-1)^{k_1} h_3 (n-u)\varphi _1\Big)\cdot
\\
&
\cdot \mathcal O\Big(\frac{1}{H}\sum\limits_{|h_4|\leq H}c_ {H}^{(j_2)} (h_4) \sum\limits_{u} e \Big( (-1)^{k_1} h_4(n-u)\varphi _1 \Big)\Big),
\nonumber
\end{split}
\end{equation*}
pour tout  $0\leq u<q_{1,k_1-1}$ ou $q_{1,k_1-1}\leq u<q_{1,k_1}$ et $i_1,i_2$ et $j\in\{1,2\}$.
En majorant trivialement sur la somme sur  $r$ puis \'{e}ventuellement  le premier terme d'erreur  par $q_{1,k_1}$.
De plus, nous utilisons l'estimation triviale $|b_ {H}^{(i)} (h_1)M^{(i)}_{k_2}(h_1,\beta)|\le 1$ et le fait que $\mathcal O$-termes sont des nombres r\'{e}els nonn\'{e}gatifs. Cel\`{a} emm\`{e}ne \`{a}
 \begin{eqnarray}
\begin{split}
S_{2,1,1}\ll&\frac{N}{H}q_{1,k_1}
 \sum\limits_{{|h_1|, |h_2|\leq H}\atop{h_1+h_2\neq0}}\sum\limits_{|h_4|\leq H}
\Big|\sum_{n\leq N}e \Big((-1)^{k_2} (h_1+h_2)n\varphi_2+(-1)^{k_1} h_4n\varphi _1\Big)\Big|
\\
&\qquad\qquad\qquad\qquad\qquad
\cdot\Big|\sum\limits_{u} e \Big( -(-1)^{k_1} h_4u\varphi _1\Big)\Big|
\\
\ll&\frac{N}{H}q_{1,k_1}
 \sum\limits_{{|h_1|, |h_2|\leq H}\atop{h_1+h_2\neq0}}\sum\limits_{|h_4|\leq H}
\Big|\sum_{n\leq N}e \Big((-1)^{k_2} (h_1+h_2)n\varphi_2+(-1)^{k_1}h_4n\varphi _1\Big)\Big|
\\
& \qquad\qquad\qquad\qquad\qquad
\cdot\min\Big(q_{1,k_1}, \Big\|h_4\varphi _1\Big\|^{-1}\Big)
 \\
\ll&\frac{N}{H} Hq_{1,k_1}^2
\sum\limits_{{|h_1|, |h_2|\leq H}\atop{h_1+h_2\neq0}}
\sup_{h_4\in\mathbb{Z}}\Big|\sum_{n\leq N} e \Big(n\Big((-1)^{k_2} (h_1+h_2)\varphi_2 +(-1)^{k_1} h_4\varphi _1\Big)\Big)\Big|
\\
\ll&Nq_{1,k_1}^2
\sum\limits_{1\leq|h|\leq2 H}(1+2H-|h|)
\sup_{h_4\in\mathbb{Z}}\min \Big(N,\Big\| h\varphi_2+h_4\varphi _1\Big\|^{-1}\Big).\nonumber
\end{split}
\end{eqnarray}
Par suite, en employant l'in\'{e}galit\'{e} suivante qui d\'{e}coule directement du  Lemme \ref{lem2.2}, de sorte que
\begin{equation}\label{11}
\Big| h\varphi_2+ h_4\varphi _1\Big|>c \max(|h|,|h_4|)^{-2-\varepsilon},
\end{equation}
nous obtenons
\begin{equation}\label{eq:S211final}
   S_{2,1,1}\ll Nq_{1,k_1}^2 \sum\limits_{1\leq|h|\leq2 H}  H H^{2+\varepsilon}
\ll N  q_{1,k_1}^2  H^{4+\varepsilon}.
\end{equation}

\medskip

{\parindent 0pt \it \ref{sec:h1h1not0}.2. Si les deux facteurs contribu\'es par $E_1$ sont un terme d'erreur et un terme principal.}
Sans restreindre la  g\'{e}n\'{e}ralit\'{e}, nous pouvons supposer que le second facteur est   un terme  d'erreur
 \begin{eqnarray}
\begin{split}
S_{2,1,2}\ll\frac{N}{R}\sum\limits_{0\leq |r|<R}\Big(1-\frac{|r|}{R}\Big)
\sum_{n\leq N} &\sum\limits_{{|h_1|, |h_2|\leq H}\atop{h_1+h_2\neq0}} b_ {H}^{(i_1)} (h_1)M^{(i_1)}_{k_2}(h_1,\beta)
e \Big((-1)^{k_2} h_1(n+r)\varphi_2\Big)\nonumber\\
&\qquad \quad\quad \cdot \overline{b_ {H}^{(j_1)} (-h_2)} \overline{M^{(j_1)}_{k_2}(-h_2,\beta)}e \Big(-(-1)^{k_2} h_2n\varphi_2\Big)
\\
& \cdot \sum\limits_{|h_3|\leq H} b_ {H}^{(i_2)} (h_3)M^{(i_2)}_{k_1}(h_3,\beta)
e \Big((-1)^{k_1} h_3(n+r)\varphi_1\Big)\nonumber\\
&\cdot \mathcal O\Big(\frac{1}{H}\sum\limits_{|h_4|\leq H}{c_ {H}^{(j_2)} (h_4)} \sum\limits_{u} e \Big( (-1)^{k_1} h_4 (n-u)\varphi _1\Big)\Big),
\nonumber
\end{split}
\end{eqnarray}
pour tout  $0\leq u<q_{1,k_1-1}$ ou $q_{1,k_1-1}\leq u<q_{1,k_1}$ et $i_1,i_2,j_1,j_2\in\{1,2\}.$
Dans cette estimation, nous majorons trivialement  la somme sur  $r$ par $R$, le terme $M^{(i_2)}_{k_1}(h_3,\beta)$ par $q_{1,k_1}$, et les termes $|b_ {H}^{(i_1)} (h_1)M^{(i_1)}_{k_2}(h_1,\beta)|$ et $|\overline{b_ {H}^{(j_1)} (-h_2)}\overline{M^{(j_1)}_{k_2}(-h_2,\beta)}|$ par $1$ et nous utilisons \eqref{10} pour $b_ {H}^{(i_2)} (h_3)$. Alors, il suit
 \begin{eqnarray*}
\begin{split}
S_{2,1,2}\ll&\frac{N}{H}q_{1,k_1}\sum\limits_{|h_3|\leq H}\min \Big(\frac{1}{|h_3|},\varphi _1^{-k_1}\Big)& \\
&\cdot\sum\limits_{{|h_4|\leq H}} \sum\limits_{{|h_1|, |h_2|\leq H}\atop{h_1+h_2\neq0}}\Big|\sum_{n\leq N}
e \Big((-1)^{k_2} (h_1+h_2)n\varphi_{2 }+(-1)^{k_1} h_4n\varphi_1\Big)\Big| \Big|\sum\limits_{u} e \Big( -(-1)^{k_1} h_4u\varphi _1\Big)\Big|.
\end{split}
\end{eqnarray*}
En majorant la somme sur $u$ trivialement  la somme sur $h_3$ par $\log H$, nous \'{e}crivons
 \begin{eqnarray}
   S_{2,1,2}&\ll&
N \log H q_{1,k_1}^2
\sum\limits_{{|h_1|, |h_2|\leq H}\atop{h_1+h_2\neq0}}
\sup_{h_4\in\mathbb{Z}}\Big|\sum_{n\leq N} e \Big(n((-1)^{k_2} (h_1+h_2)\varphi_2+(-1)^{k_1} h_4\varphi _1)\Big)
\Big|\nonumber\\
&\ll&N\log H q_{1,k_1}^2
\sum\limits_{1\leq|h|\leq2 H}(1+2H-|h|)
\sup_{h_4\in\mathbb{Z}}\min \Big(N,\Big\|h\varphi_2+ h_4\varphi _1\Big\|^{-1}\Big).\nonumber
\end{eqnarray}
Par cons\'{e}quent, un appel \`{a} (\ref{11}) conduit \`{a}
 \begin{equation}\label{eq:S212final}
   S_{2,1,2}\ll N\log Hq_{1,k_1}^2  \sum\limits_{|h|\leq2 H} H H^{2+\varepsilon}
\ll N q_{1,k_1}^2H^{4+\varepsilon} \log H .
\end{equation}

\medskip

{\parindent 0pt \it \ref{sec:h1h1not0}.3. Si les deux facteurs contribu\'es par $E_1$ sont des termes principaux.}
Dans ce cas, nous commen\c cons par
\begin{equation*}
\begin{split}
S_{2,1,3}=\frac{N}{R}\sum\limits_{0\leq |r|<R}\Big(1-\frac{|r|}{R}\Big)
\sum_{n\leq N} &\sum\limits_{|h_3|\leq H} b_ {H}^{(i_1)} (h_3)M^{(i_1)}_{k_1}(h_3,\vartheta)
e \Big((-1)^{k_1} h_3(n+r)\varphi_1\Big)\nonumber\\
\cdot&\sum\limits_{|h_4|\leq H} \overline{b_ {H}^{(j_1)} (-h_4)} \overline{M^{(j_1)}_{k_1}(-h_4,\vartheta)}
e \Big((-1)^{k_1} h_4n\varphi_{1 }\Big)\nonumber\\
\cdot&\sum\limits_{|h_1|\leq H} b_ {H}^{(i_2)} (h_1)M^{(i_2)}_{k_2}(h_1,\vartheta)
e \Big((-1)^{k_2} h_1(n+r)\varphi_2\Big)
\\
\cdot&\sum\limits_{|h_2|\leq H} \overline{b_ {H}^{(j_2)} (-h_2)} \overline{M^{(j_2)}_{k_2}(-h_2,\beta)}
e \Big((-1)^{k_2} h_2n\varphi_2\Big),
\end{split}
\end{equation*}
que nous pouvons r\'{e}\'{e}crire comme \'{e}tant
\begin{equation*}
\begin{split}
&S_{2,1,3}=\frac{N}{R}\sum\limits_{0\leq |r|<R}\Big(1-\frac{|r|}{R}\Big)
\sum_{n\leq N}\nonumber\\
&\quad\sum\limits_{{|h_3|, |h_4|\leq H}}b_ {H}^{(i_1)} (h_3)M^{(i_1)}_{k_1}(h_3,\vartheta)\overline{b_ {H}^{(j_1)} (-h_4)} \overline{M^{(j_1)}_{k_1}(-h_4,\beta)}e \Big((-1)^{k_1} (h_3+h_4)n\varphi_1+(-1)^{k_1}h_3r\varphi_1\Big)\nonumber\\
&\quad\sum\limits_{{|h_1|, |h_2|\leq H}\atop{h_1+h_2\neq0}}
b_ {H}^{(i_2)} (h_1)M^{(i_2)}_{k_2}(h_1,\beta)\overline{b_{H}^{(j_2)} (-h_2)} \overline{M^{(j_2)}_{k_2}(-h_2,\beta)}
e \Big((-1)^{k_2} (h_1+h_2)n\varphi_2+(-1)^{k_2} h_1r\varphi_2\Big).
\nonumber
\end{split}
\end{equation*}
Nous estimons trivialement  la somme sur $r$ ainsi que tous les termes   $M_{k_1} $ et $M_{k_2} $   et nous appelons l'in\'{e}galit\'{e} (\ref{02}) pour tous les termes $b_ {H}$,  ce qui nous emm\`{e}ne  \`{a}
  \begin{eqnarray}
S_{2,1,3}&\ll&N\sum\limits_{{|h_1|, |h_2|\leq H}\atop{h_1+h_2\neq0}}\sum\limits_{{|h_3|, |h_4|\leq H}}\Big|
\sum_{n\leq N}e \Big((-1)^{k_1}(h_3+h_4)n\varphi_1+(-1)^{k_2} (h_1+h_2)n\varphi_2\Big)\Big|\nonumber\\
&\ll&N\sum\limits_{{1\leq|h|\leq 2H \atop |l|\leq 2H}}(1+2H-|h|)(1+2H-|l|)
\Big|\sum_{n\leq N}e \Big(h n\varphi_1+ l n\varphi_2\Big)\Big|.\nonumber
\end{eqnarray}
\`{A} pr\'{e}sent, nous appliquons le  Lemme~\ref{lem2.2} afin d'estimer  la somme sur  $n$ et nous obtenons
\begin{equation}\label{eq:S213final}
\begin{split}
S_{2,1,3}\ll&N\sum\limits_{{1\leq|h|\leq 2H \atop |l|\leq 2H}}(1+2H-|h|)(1+2H-|l|)
\min\Big(N,\Big\|h\varphi_1+l\varphi_2\Big\|^{-1}\Big)\\
\ll&N H^2\sum\limits_{{1\leq|h|\leq 2H \atop |l|\leq 2H}}\max(|h|, |l|)^{2+\varepsilon}
\\
\ll&N H^{6+\varepsilon}.
\end{split}
\end{equation}
En combinant \eqref{eq:S211final}, \eqref{eq:S212final} et \eqref{eq:S213final} , nous arrivons \`{a}
\begin{equation}\label{eq:S21final}
S_{2,1} \ll
Nq_{1,k_1}^2 H^{4+\varepsilon}
+ Nq_{1,k_1}^2H^{4+\varepsilon} \log H
+ NH^{6+\varepsilon} \ll Nq_{1,k_1}^2H^{7}.
\end{equation}
Nous finissons  le  Cas 2 avec l'estimation suivante, en mettant la derni\`{e}re majoration avec  \eqref{eq:S20final},
\begin{equation}\label{eq:S2final}
\begin{split}
S_2 &\ll N^2 e^{-c_2k_2}\log H\Big(\frac{q_{2,k_2}}{R} + \log q_{2,k_2} \Big) +NHq_{1,k_1} + Nq_{1,k_1}^2H^{7} \\
&\ll N^2 e^{-c_2k_2}\log H\Big(\frac{q_{2,k_2}}{R} + \log q_{2,k_2} \Big)  + Nq_{1,k_1}^2H^{7}.
\end{split}
\end{equation}
\subsection{Cas 3 :
Si un des facteurs contribu\'es par $E_2$ est un terme  d'erreur et l'autre est un  terme principal.}
Sans perte de g\'{e}n\'{e}ralit\'{e}, nous assumons que le second facteur est   un terme  d'erreur.
Alors,
\begin{eqnarray*}
\begin{split}
S_3=\frac{N}{R}\sum\limits_{0\leq |r|<R}\Big(1-\frac{|r|}{R}\Big)
&\sum_{n\leq N} e \Big(\vartheta \Big(S_{\alpha_1,k_1}(n+r)- S_{\alpha_1,k_1}(n)\Big)\Big)\\
\cdot&\sum\limits_{|h_1|\leq H} b_ {H}^{(i_1)} (h_1)M^{(i_1)}_{k_2}(h_1,\beta)
e \Big((-1)^{k_2} h_1(n+r)\varphi_2\Big)
\\
\cdot& \,\mathcal O\Big(\frac{1}{H}\sum\limits_{|h_2|\leq H}c_ {H}^{(j_1)} (h_2) \sum\limits_{u} e \Big( (-1)^{k_2} h_2 (n-u)\varphi _2\Big)\Big),
\end{split}
\end{eqnarray*}
pour tout  $0\leq u<q_{2,k_2-1}$ ou $q_{2,k_2-1}\leq u<q_{2,k_2}$ et $i_1,j_1\in\{1,2\}$ et
l' expression dans le terme d'erreur est un nombre r\'{e}el  non n\'{e}gatif. En suivant le  cas pr\'{e}c\'{e}dant et en rempla\c{c}ant les facteurs  d'exponentielles dans
$E_1$ par leurs expressions dans la Proposition~\ref{prop1}, nous obtenons encore une fois pour ce sous-cas seize sommes de produits des  termes principaux et des termes d'erreurs dont nous distinguons  aussi pour $S_3$ trois cas correspondants aux nombres des termes d'erreur contribu\'{e}s par $E_1$. Nous utilisons la  d\'{e}composition suivante
\[
S_{3}=S_{3,1}+S_{3,2}+S_{3,3}.
\]
Dor\'{e}navant, chaque sommation sera trait\'{e}e d'une fa\c con s\'{e}par\'{e}e
\subsubsection{Si les deux facteurs contribu\'es par $E_1$ sont des termes d'erreurs}
Das ce cas, la contribution de  $S_{3,1}$ \`{a} $S_3$ est
\begin{eqnarray*}
\begin{split}
S_{3,1}=\frac{N}{R}\sum\limits_{0\leq |r|<R}\Big(1-\frac{|r|}{R}\Big)&\sum_{n\leq N}\mathcal O\Big(\frac{1}{H}\sum\limits_{|h_3|\leq H}c_ {H}^{(j_1)} (h_3) \sum\limits_{u} e \Big( (-1)^{k_1} h_3 (n-u)\varphi _1\Big)\Big)\\
\cdot & \,\mathcal O\Big(\frac{1}{H}\sum\limits_{|h_4|\leq H}c_ {H}^{(j_2)} (h_4) \sum\limits_{u} e \Big( (-1)^{k_1} h_4(n-u)\varphi _1\Big)\Big)\\
\cdot &\sum\limits_{|h_1|\leq H} b_ {H}^{(i_1)} (h_1)M^{(i_1)}_{k_2}(h_1,\beta)
e \Big((-1)^{k_2} h_1(n+r)\varphi_2\Big)
\\
\cdot & \,\mathcal O\Big(\frac{1}{H}\sum\limits_{|h_2|\leq H}c_ {H}^{(j_3)} (h_2) \sum\limits_{u} e \Big( (-1)^{k_2} h_2 (n-u)\varphi _2\Big)\Big),\nonumber
\end{split}
\end{eqnarray*}
pour tout $0\leq u<q_{1,k_1-1}$ ou $q_{1,k_1-1}\leq u<q_{1,k_1}$ et $0\leq u<q_{2,k_2-1}$ ou $q_{2,k_2-1}\leq u<q_{2,k_2}$ et $i_1,j_1,j_2,j_3\in\{1,2\}.$
Nous majorons trivialement la somme sur  $r$, aussi le premier terme d'erreur par $q_{1,k_1}$ et nous appelons
$|b_ {H}^{(i_1)}(h_1) M^{(i_1)}_{k_2}(h_1,\beta)|\leq 1$ ce qui nous conduit \`{a}
\begin{eqnarray*}
\begin{split}
S_{3,1}\ll&\frac{N}{H^2}\sum\limits_{|h_1|\leq H} q_{1,k_1}
\sum\limits_{|h_2|,|h_4|\leq H}\Big|\sum_{n\leq N}e \Big((-1)^{k_2} h_2 n\varphi _2+ (-1)^{k_1} h_4n\varphi _1\Big)\Big|\nonumber\\
& \qquad\cdot \Big|\sum\limits_{u_1} e \Big( (-1)^{k_1} h_4u_1\varphi _1\Big)\Big|
\Big|\sum\limits_{u_2} e \Big( (-1)^{k_2} h_2 u_2\varphi _2\Big)\Big|,\\
\ll&\frac{N}{H}q_{1,k_1}^2 q_{2,k_2}\sum\limits_{|h_2|,|h_4|\leq H}\Big|\sum_{n\leq N}e \Big((-1)^{k_2} h_2 n\varphi _2+ (-1)^{k_1} h_4n\varphi _1\Big)\Big|\\
\ll&\frac{N}{H}q_{1,k_1}^2 q_{2,k_2}\bigg(N+ \sum\limits_{|h_2|,|h_4|\leq H \atop (h_2,h_4)\not=0}\Big|\sum_{n\leq N}e \Big((-1)^{k_2} h_2 n\varphi _2+ (-1)^{k_1} h_4n\varphi _1\Big)\Big| \bigg)\\
\ll&\frac{N}{H} q_{1,k_1}^2 q_{2,k_2}\bigg(N+ \sum\limits_{|h_2|,|h_4|\leq H \atop (h_2,h_4)\not=0} \min\Big(N,\| h_2\varphi _2+ h_4\varphi _1\|^{-1}\Big) \bigg).
\end{split}
\end{eqnarray*}
Si  $(h_2,h_4)\not=0$, en tenant compte  de  l'in\'{e}galit\'{e} (\ref{11}),  il s'ensuit
\begin{eqnarray}\label{eq:S31final}
\begin{split}
S_{3,1}&\ll \frac{N}{H}q_{1,k_1} q_{2,k_2}
\bigg(N + \sum\limits_{|h_2|,|h_4|\leq H \atop (h_2,h_4)\not=0} \min\Big(N,\| h_2\varphi _2+ h_4\varphi _1\|^{-1}\Big)
\bigg) \\
&\ll
\frac{N}{H} q_{1,k_1}^2 q_{2,k_2}\bigg(N +
\sum\limits_{|h_2|,|h_4|\leq H} \max(|h_2|,|h_4|)^{2+\varepsilon}\bigg)  \\
&\ll  Nq_{1,k_1}^2 q_{2,k_2} \Big( \frac{N}{H} + H^{3+\varepsilon}  \Big).
\end{split}
\end{eqnarray}

\subsubsection{Si les deux facteurs contribu\'es par $E_1$ sont un terme d'erreur et un terme principal}
Sans perte de g\'{e}n\'{e}ralit\'{e}, nous supposons que le second facteur est un terme  d'erreur.
La sommation correspondante \`{a} $S_3$ est alors lue comme suit
\begin{equation*}
\begin{split}
S_{3,2}=\frac{N}{R}\sum\limits_{0\leq |r|<R}\Big(1-\frac{|r|}{R}\Big)&
\sum_{n\leq N}\sum\limits_{|h_1|\leq H} b_ {H}^{(i_1)} (h_1)M^{(i_1)}_{k_2}(h_1,\beta)
e \Big((-1)^{k_2} h_1(n+r)\varphi_2\Big)\\
\cdot&\,\mathcal O\Big(\frac{1}{H}\sum\limits_{|h_2|\leq H}c_ {H}^{(j_1)} (h_2) \sum\limits_{u_1} e \Big( (-1)^{k_2} h_2 (n-u_1)\varphi _2\Big)\Big)\\
\cdot&\sum\limits_{|h_3|\leq H} b_ {H}^{(i_2)} (h_3)M^{(i_2)}_{k_1}(h_3,\beta)
e \Big((-1)^{k_1} h_3(n+r)\varphi_1\Big) \\
\cdot&\,\mathcal O\Big(\frac{1}{H}\sum\limits_{|h_4|\leq H}c_ {H}^{(j_2)} (h_4) \sum\limits_{u} e \Big( (-1)^{k_1} h_4 (n-u_2)\varphi _1\Big)\Big),
\end{split}
\end{equation*}
pour tout $0\leq u<q_{1,k_1-1}$ ou $q_{1,k_1-1}\leq u<q_{1,k_1}$ et   $0\leq u<q_{2,k_2-1}$ ou $q_{2,k_2-1}\leq u<q_{2,k_2}$ et $i_1,j_1,i_2,j_2\in\{1,2\}$.\\
Maintenant, il suffit d'estimer  la somme  sur $r$ trivialement  ainsi que les termes $M^{(i_2)}_{k_1}$ et $M^{(i_1)}_{k_2}$. Gr\^{a}ce \`{a} la propri\'{e}t\'{e} (\ref{10}), nous majorons  les  $b_ {H}$ par $\frac{1}{|h|}$ et nous constatons
 \begin{eqnarray*}
\begin{split}
S_{3,2}\ll&\frac{N}{H^2}\sum\limits_{|h_3|\leq H} \frac{1}{|h_3|} q_{1,k_1}
\sum\limits_{|h_1|\leq H} \frac{1}{|h_1|}  q_{2,k_2}
\sum\limits_{|h_2|,|h_4|\leq H}
\Big|\sum_{n\leq N} e \Big((-1)^{k_2} h_2n\varphi_2+(-1)^{k_1} h_4n\varphi_1\Big)\Big|\\
&\quad\cdot\,\Big|\sum\limits_{u_1} e \Big( -(-1)^{k_2}  h_2u_1\varphi _2\Big)\Big|
\Big|\sum\limits_{u_2} e \Big( -(-1)^{k_1}  h_4u_2\varphi _1\Big)\Big|\\
\ll&\frac{N}{H^2}(\log H )^2 q_{1,k_1}^2q_{2,k_2}^2\bigg( N +
\sum\limits_{|h_2|,|h_4|\leq H \atop (h_2,h_4) \not=0}
\min \Big(N,\|h_2\varphi_2+h_4\varphi_1\|^{-1}\Big)\bigg).
\end{split}
\end{eqnarray*}
Encore une fois, nous employons (\ref{11}) qui d\'{e}coule directement de l'in\'{e}galit\'{e} de Schmidt. Donc, nous parvenons \`{a} \'{e}crire
 \begin{equation}\label{eq:S32final}
 \begin{split}
   S_{3,2}&\ll\frac{N}{H^2}(\log H )^2 q_{1,k_1}^2q_{2,k_2}^2\bigg(N+\sum\limits_{|h_2|,|h_4|\leq H \atop (h_2,h_4) \not=0} \max(|h_2|,|h_4|)^{2+\varepsilon}\bigg) \\
  & \ll N(\log H )^2 q_{1,k_1}^2q_{2,k_2}^2 \Big(\frac{N}{H^2}+H^{2+\varepsilon}\Big).
\end{split}
\end{equation}

\subsubsection{Si les deux facteurs contribu\'es par $E_1$ sont des termes principaux}\label{sec:S33}
Dans ce cas, la sommation correspondante \`{a}
 $S_3$ est alors lue comme suit
\begin{eqnarray*}
\begin{split}
S_{3,3}=\frac{N}{R}\sum\limits_{0\leq |r|<R}\Big(1-\frac{|r|}{R}\Big)&
\sum_{n\leq N} \sum\limits_{|h_3|\leq H} b_ {H}^{(i_1)} (h_3)M^{(i_1)}_{k_1}(h_3,\vartheta)
e \Big((-1)^{k_1} h_3(n+r)\varphi_1\Big)\nonumber\\
\cdot&\sum\limits_{|h_4|\leq H} \overline{b_ {H}^{(j_1)} (-h_4)} \overline{M^{(j_1)}_{k_1}(-h_4,\vartheta)}
e \Big((-1)^{k_1} h_4n\varphi_1\Big)\\
\cdot&\sum\limits_{|h_1|\leq H} b_ {H}^{(i_2)} (h_1)M^{(i_2)}_{k_2}(h_1,\beta)
e \Big((-1)^{k_2} h_1(n+r)\varphi_2\Big)
\\
\cdot&\,\mathcal O\Big(\frac{1}{H}\sum\limits_{|h_2|\leq H}c_ {H}^{(j_2)} (h_2) \sum\limits_{u} e \Big( (-1)^{k_2} h_2 (n-u)\varphi _2\Big)\Big),
\end{split}
\end{eqnarray*}
pour tout  $0\leq u<q_{2,k_2-1}$ ou $q_{2,k_2-1}\leq u<q_{2,k_2}$ et $i_1,j_1,i_2,j_2\in\{1,2\}$.
Pour traiter ce cas, nous majorons  trivialement la somme sur $r$ ainsi que les termes $M^{(i)}_{k}$ et $|b_ {H}^{(j_1)} M^{(j_1)}_{k_1}|\leq 1$. Par cons\'{e}quent, nous  obtenons en utilisant  \eqref{10} pour  $b_ {H}^{(i_1)} (h_3)$,  $b_ {H}^{(i_2)} (h_1)$ et le Lemme~\ref{lemsch}
\begin{equation}\label{eq:S33final}
\begin{split}
S_{3,3}&\ll \frac{N}{H}\sum\limits_{{|h_3|\leq H}} \frac{1}{|h_3|}q_{1,k_1}
\sum\limits_{|h_1|\leq H} \frac{1}{|h_1|}q_{2,k_2}  \\
&\sum\limits_{0\leq|h_2|,|h_4|\leq H}\Big|\sum\limits_{n\leq N} e \Big( (-1)^{k_2} h_2n\varphi _2+(-1)^{k_1} h_4n\varphi _1\Big)\Big|
\cdot\Big|\sum\limits_{u} e \Big( -(-1)^{k_2}  h_2u\varphi _2\Big)\Big| \\
&\ll\frac{N\log^2 H}{H}q_{1,k_1}q_{2,k_2}^2 \bigg(N+ \sum\limits_{0\le |h_2|,|h_4|\leq H\atop{(h_2,h_4)\neq 0}} \min \Big(N,\|h_2\varphi_2+h_4\varphi_1\|^{-1}\Big)
 \bigg)\\
&\ll N(\log^2 H)q_{1,k_1}q_{2,k_2}^2 \Big(\frac{N}H + H^{3+\varepsilon}\Big).
\end{split}
\end{equation}

L'assemblement de  \eqref{eq:S31final}, \eqref{eq:S32final} et \eqref{eq:S33final} am\`{e}ne \`{a}
\begin{equation}\label{eq:S3final}
\begin{split}
S_3 \ll &Nq_{1,k_1}^2 q_{2,k_2} \Big( \frac{N}{H} + H^{3+\varepsilon}  \Big) + N(\log H )^2 q_{1,k_1}^2q_{2,k_2}^2 \Big(\frac{N}{H^2}+H^{2+\varepsilon}\Big) \\
&+ N(\log^2 H)q_{1,k_1}q_{2,k_2}^2 \Big(\frac{N}H + H^{3+\varepsilon}\Big) \\
\ll& Nq_{1,k_1}^2q_{2,k_2}^2\log^2H \Big(\frac{N}{H}+ H^{3+\varepsilon} \Big)
\ll Nq_{1,k_1}^2q_{2,k_2}^2\log^2H \Big(\frac{N}{H}+ H^{4} \Big)
\end{split}
\end{equation}
ce qui conclut ce cas.

\subsection{Conclusion de la preuve du Th\'{e}or\`{e}me~\ref{th}}
Finallement, en  combinant les trois cas \eqref{eq:S1final}, \eqref{eq:S2final}, \eqref{eq:S3final} et  en  mettant  les sept termes d'erreur obtenus dans l'in\'{e}galit\'{e} (\ref{eqp}) et en rempla\c{c}ant chaque terme par sa valeur, nous aboutissons   \`{a}
\begin{eqnarray*}
\begin{split}
   \Big|\sum_{ n<N}e (\vartheta S_{\alpha_1}(n)+\beta S_{\alpha_2}(n))\Big|^2
\ll& \;
N^2 q^{-1/2}_{1,k_1}q^{3/2}_{2,k_2}+N^2q^{2}_{2,k_2}\frac{\log H}{H}+ Nq_{1,k_1}q_{2,k_2}^2
\\
&+N^2 e^{-c_2k_2}\log H\Big(\frac{q_{2,k_2}}{R} + \log q_{2,k_2} \Big)  + Nq_{1,k_1}^2H^{7}
\\
&+Nq_{1,k_1}^2q_{2,k_2}^2\log^2H \Big(\frac{N}{H}+ H^{4} \Big)
\\
&+NR+\frac{N^2R}{q_{1,k_1-1}}+\frac{N^2R}{q_{2,k_2-1}}.
\end{split}
\end{eqnarray*}
Par suite, nous effectuons le  choix  d'une nouvelle variable $\sigma$  assez petite pour
$$
k_1=\Big\lfloor 4\sigma \frac{\log N}{\log \varphi_1}\Big\rfloor ,\
k_2=\Big\lfloor \frac{\sigma\log N}{\log \varphi_{_2}}\Big\rfloor, \
H=\Big\lfloor N^{16\sigma-c_2 \sigma/\log\varphi_{2} }\Big\rfloor, \
R=\Big\lfloor N^{\sigma-c_2 \sigma/2\log\varphi_{_2} } \Big \rfloor.
$$
Ces choix conduisent \`{a}
\begin{eqnarray*}
\begin{split}
  \Big|\sum_{ n<N}e (\vartheta S_{\alpha_1}(n)+\beta S_{\alpha_2}(n))\Big|^2
\ll& N^{2-\sigma/2 }+\log N N^{2-14\sigma+c_2 \sigma/\log\varphi_{2}}+N^{1+6\sigma}\\
+&\log N N^{2-c_2 \sigma/2\log\varphi_{2}}+\log^2 N N^{2-c_2 \sigma/\log\varphi_{2}}   \\
&\qquad+N^{116\sigma-7c_2 \sigma/\log\varphi_{2}}\\
+&\log^2 N (N^{2-6\sigma +c_2 \sigma/\log\varphi_{2}}+N^{1+58\sigma -4c_2 \sigma/\log\varphi_{2}})\\
&\qquad+N^{1+\sigma-c_2 \sigma/2\log\varphi_{2}}\\
+&N^{2-3\sigma-c_2 \sigma/2\log\varphi_{2}}+N^{2-c_2 \sigma/2\log\varphi_{2}} \\
\ll&N^{2(1-\delta)},
\end{split}
\end{eqnarray*}
pour certaine  $\delta > 0$ . Ce qui  ach\`eve \`{a} la contribution de l'estimation de la somme d'exponentielles d\'{e}sir\'{e}e et le Th\'{e}or\`{e}me~\ref{th} est bien prouv\'{e}.

\subsection{Preuve du  Corollaire~\ref{cor}}
Nous utilisons la relation d'orthogonalit\'{e} classique \`{a} savoir
\begin{equation*}
\frac{1}{b}\sum_{j=0}^{b-1}
e\Big(\frac{j(a-c)}{b}\Big)=
 \left\{
\begin{array}{lll}
1,  &si&c\equiv a\pmod{b}\\
0 ,  &sinon&
\end{array} \right.
\end{equation*}
En isolant le terme correspondant \`{a} $j_1=j_2=0$ et en vue du Th\'{e}or\`{e}me \ref{th}, nous obtenons
\begin{align*}
\Big| \{ 0\leq n<N;\ &S_{\alpha_1}(n)\equiv a_1\pmod{b_1},\ S_{\alpha_2}(n)\equiv a_2\pmod{b_2}\}\Big|\\
 =&\frac{1}{b_1b_2}\sum_{n<N}\sum_{0\leq j_1<b_1}e\Big(\frac{j_1(S_{\alpha_1}(n)-a_1)}{b_1}\Big)
\sum_{0\leq j_2<b_2}e\Big(\frac{j_2(S_{\alpha_2}(n)-a_2)}{b_2}\Big)\nonumber\\
 =&\frac{1}{b_1b_2}\sum_{{0\leq j_1<b_1}\atop{0\leq j_2<b_2}}e\Big(-\frac{j_1a_1}{b_1}
-\frac{j_2a_2}{b_2}\Big) \sum_{n<N}e\Big(\frac{j_1S_{\alpha_1}(n)}{b_1}+\frac{j_2S_{\alpha_2}(n)}{b_2}\Big) \nonumber\\
=&\frac{N}{b_1b_2}+\frac{1}{b_1b_2}
\sum_{{0\leq j_1<b_1}\atop{{0\leq j_2<b_2}\atop{(j_1,j_2)\neq(0,0)}}}e\Big(-\frac{j_1a_1}{b_1}
-\frac{j_2a_2}{b_2}\Big) \sum_{n<N}e\Big(\frac{j_1S_{\alpha_1}(n)}{b_1}+\frac{j_2S_{\alpha_2}(n)}{b_2}\Big)
   \nonumber\\
=&\frac{N}{b_1b_2}
+\mathcal O\Big(\frac{1}{b_1b_2}  \sum_{{0\leq j_1,j_2<b_1}\atop{(j_1,j_2)\not=0}}
\Big|\sum_{n<N}e\Big(\frac{j_1S_{\alpha_1}(n)}{b_1}+\frac{j_2S_{\alpha_2}(n)}{b_2}\Big) \Big|  \Big)  \\
=&
\frac{N}{b_1b_2}+\mathcal O(N^{1-\delta}).
\end{align*}
Il est \`{a} signaler que les conditions $(b_1,m_1)=1$ ou  $(b_2,m_2)=1$ sont n\'{e}cessaires pour \'{e}viter la possibilit\'{e} $\frac{j_i m_i}{b_i}$ \'{e}tant un entier $(i=1,2)$,  ce qui contredit la condition donn\'{e}e dans le Th\'{e}or\`{e}me \ref{th}.
D'o\`{u} le  corollaire est achev\'{e}.\\

\bibliographystyle{plain}
\bibliography{Bib5}

\end{document}